\documentclass[a4paper,11pt,dvips]{article}

\usepackage{amsmath,amsfonts,amsthm,amssymb}
\usepackage[T1]{fontenc}
\usepackage{mathptmx}
\usepackage{mhequ}
\usepackage{mhsymb}
\usepackage{mathptmx}
\usepackage{cite}
\usepackage{bm} 						
                						
\usepackage[colorlinks=true,
						linkcolor=black,
						citecolor=black,
						anchorcolor=black,
						urlcolor=black,
						filecolor=black,
						menucolor=black,
						runcolor=black]{hyperref}	 

\usepackage[subnum]{cases}
\usepackage{subfigure}
\usepackage{appendix}
\usepackage{pdftricks}
\usepackage{pst-node}
\usepackage[vmargin=3.5cm,hmargin=3.5cm]{geometry}  
\usepackage{microtype}
	\DisableLigatures{encoding = *, family = * } 		

\newcommand{\I}{\mathbb{I}} 
\newcommand{\SO}[1]{\operatorname{SO}(#1)}
\newcommand{\Real}[1]{\mathbb{R}^{#1}}
\newcommand{\cay}{\text{cay}}
\newcommand{\dist}{\textrm{dist}}
\newcommand{\trace}{\textrm{tr}}
\DeclareMathOperator{\D}{\textrm{D}} 	

\global \def \figscale{.48}

\title{A counterexample showing the semi-explicit
Lie-Newmark algorithm is not variational}

\author{Nawaf Bou-Rabee\thanks{Department of Mathematics, 
															Courant Institute of Mathematical Sciences, 
															New York university. E-mail: \texttt{nawaf@cims.nyu.edu}.} 
				\and
				Giulia Ortolan\thanks{Department of Information Engineering, University of Padova. 
															E-mail: \texttt{ortolang@dei.unipd.it}.}
				\and
				Alessandro Saccon\thanks{Instituto de Sistemas e Robótica, Instituto Superior Técnico, 
																	Universidade Técnica de Lisboa.
																	E-mail: \texttt{asaccon@isr.ist.utl.pt}.}}
																	
\date{}

\begin{document}
\maketitle

\begin{abstract}
This paper presents a counterexample to the conjecture that
the semi-explicit Lie-Newmark algorithm is variational.
As a consequence the Lie-Newmark method is not well-suited for long-time
simulation of rigid body-type mechanical systems.
The counterexample consists of a single rigid body in a static potential field,
and can serve as a test of the variational nature of other rigid-body integrators.

\vspace{.2cm}

\end{abstract}

%
%

\section{Introduction}

In this paper we will focus on the dynamics of a rigid body in a static potential field.
To describe this system, denote by $Q(t) \in \SO{3}$,
$W(t) = {\left[ W_1(t)\ W_2(t)\ W_3(t) \right]}^T \in \Real{3}$, and
$\I = \text{diag}(I_1,I_2,I_3)\in \Real{3\times 3}$
the configuration, body angular velocity and inertia matrix of the body, respectively.
Let $\tau: \SO{3} \to \mathbb{R}^3$ be the torque acting on the body and
$\ \widehat\,:\,\Real{3} \to \Real{3 \times 3}$ be the hat map \[
\widehat{W} = \begin{bmatrix} 0 & -W_3 & W_2 \\
             				 W_3 & 0 & -W_1 \\
					-W_2 &	W_1 & 0 \end{bmatrix}\!.
\]  In terms of this notation, the governing equations are
\begin{numcases}{\label{eq:RigidBody}}
\dot{Q} \!\!\!\!\!\!\!\!\!& $= Q \widehat{W} $\label{eq:reconstruction} \\
\I\,\dot{W} \!\!\!\!\!\!\!\!\!& $ = \I\,W \times W + \tau(Q),$  \label{eq:balance}
\end{numcases}
with initial conditions $Q(0) = Q_0 \in \SO{3}$ and $W(0) = W_0 \in \mathbb{R}^3$.
We assume this rigid body is derivable from a Lagrangian $L: \SO{3} \times \Real{3} \to \Real{}$
of the form
\begin{equation} \label{eq:Lagrangian}
L(Q,W) = T(W) - U(Q) \;,
\end{equation}
where $T(W)=\frac{1}{2}W^T \I\,W$ and $U(Q)$ are the kinetic and potential energy
of the body, respectively. This assumption implies that the torque in \eqref{eq:balance}
can be computed from the directional derivative of $U$ 
at $Q$ in the direction $Q \widehat{y}$: \[
\tau(Q)^T y =  - \D U(Q) \cdot Q\,\hat{y} \;,~~~y \in \Real{3} \;.
\]
Notice that the total energy is separable and $T(-W) = T(W)$.
The exact flow of \eqref{eq:RigidBody} possesses certain structure such as
total energy preservation, time-symmetry, and symplecticity.
Moreover, the path $Q$ lies on a configuration manifold $\SO{3}$ which possesses
a Lie-group structure.

This paper investigates the long-run behavior of two integrators
for \eqref{eq:RigidBody}: the \linebreak[4] Lie-Newmark \cite{SiVu1988, SiWo1991}
and Lie-Verlet methods \cite{BoMa2008}.
Both methods are semi-explicit, second-order accurate and symmetric.
They are also `Lie group methods' because they respect the Lie group structure of the
configuration manifold \cite{IMNZ99}.   Moreover, the complexity and implementation
of the two methods is quite similar.  The main difference between
the integrators is that the Lie-Verlet method is designed to be variational, whereas
the Lie-Newmark method is not.

Variational integrators are time-integrators adapted to the structure
of mechanical systems \cite{MaWe2001,LeMaOrWe2004a,LeMaOrWe2004b}.
They are symplectic, and in the presence of symmetry, momentum preserving.
The theory of variational integrators includes discrete analogues of Hamilton's principle,
Noether's theorem, the Euler-Lagrange equations, and the Legendre transform.
The variational nature of Lie-Verlet guarantees its excellent long-time behavior.
In fact, one can prove this.  The basic idea of the proof is to show that a trajectory of a
variational integrator is interpolated by a level set of a `modified' energy function nearby the
true energy \cite{BeGi1994, Re1999, HaLuWa2006}.   This implies that a trajectory of the
variational integrator is confined to these level sets for the duration of the simulation.
As a consequence variational integrators nearly preserve the true energy and exhibit
linear growth in global error.  For these reasons variational integrators are well-suited
for long-time simulation.

Even though the Lie-Newmark integrator is not designed to be variational,
this does not rule out the possibility that the algorithm is variational in a subtle way
like classical Newmark.   The classical Newmark family of algorithms
are widely used integrators in computational  mechanics (for an expository treatment
see, e.g., Chapter 9 of \cite{Hu1987}).    They were first proposed in \cite{Ne1959}, but
it took more than forty years for their variational nature to be established in \cite{KaMaOrWe2000}.
Specifically, Kane et al.~proved that a trajectory of the classical Newmark method is
shadowed by a trajectory of a variational algorithm.  In other words the Newmark integrators
are not symplectic, but conjugate symplectic \cite{HaLuWa2006}.
The possibility that Lie-Newmark could be analogously variational (or conjugate symplectic)
was supported by recent numerical evidence showing that the Lie-Newmark algorithm exhibits
excellent behavior akin to classical Newmark \cite{KrEn2005}.  Based on these experiments,
Krysl and Endres conjectured that the Lie-Newmark algorithm is variational.

This paper disproves this conjecture.  In particular, the paper presents a simple
numerical counterexample showing that the Lie-Newmark method exhibits systematic energy drift.
In contrast, the Lie-Verlet method nearly preserves the true energy and exhibits the qualitative
properties one expects of a variational integrator.  In summary, the Lie-Verlet method is well-suited
for long-time simulation of rigid body-type mechanical systems, while the Lie-Newmark method
is not.

The paper is organized as follows.
The algorithms used in the paper are stated in \S 2 and the
numerical `stress test' carried out in \S 3.   In \S 4 we
provide concluding remarks including other applications
of this numerical `stress test.'  Along with this paper,
we have released a simple \textsc{Matlab} implementation 
of the presented algorithms.  This release can
be retrieved from the \textsc{Matlab} Central File Exchange.

\section{Integrators}

\paragraph{Lie-Newmark}

The Lie-Newmark family of integrators was proposed more than \linebreak[4]twenty years
ago in \cite{SiVu1988}.  These methods consist of a Newmark-style discretization
of \eqref{eq:balance} and a discretization of \eqref{eq:reconstruction} that ensures
the configuration update remains on $\SO{3}$.
These methods were motivated by the need to develop conserving
algorithms that can efficiently simulate the structural dynamics of rods and shells.   For example,
consider simulating large deformations of a three-dimensional finite-strain rod model.
The rod is typically discretized using $N$ copies of $\mathbb{R}^3 \times \SO{3}$
where $N$ is the number of discretization points along the line of centroids of the rod.  The
configuration of the rod at each point along the line of centroids is specified by an element
of $\SO{3}$.  The dynamical behavior of the rod can then be estimated
by simulating the dynamics of $N$ rigid bodies with torques due to
elastic coupling between bodies.

This paper focuses on a specific member of the Lie-Newmark family
tested in \cite{KrEn2005}.  This member is the Lie group analog of the so-called explicit
Newmark method on vector spaces (see Chapter 9 of \cite{Hu1987}).  In molecular
dynamics the explicit Newmark method is known as the Verlet integrator \cite{HaLuWa2006}.
Given $(Q_k, W_k) \in \Real{3} \times \SO{3}$ and time-stepsize $h$, the Lie-Newmark algorithm
determines  $(Q_{k+1}, W_{k+1})$ by the following iteration rule:
\begin{numcases}{\label{eq:LieNewmark}}
	W_{k+\frac{1}{2}} =
		W_k + \frac{h}{2} \I^{-1} \left( \I\,W_k \times W_k + \tau(Q_k)\right)  \;,	\label{eq:LN_1}	 & \\
	Q_{k+1} =
		Q_k\,\cay(h W_{k+\frac{1}{2}})  \;,	\label{eq:LN_2} 	&  \\
	W_{k+1} =
		W_{k+\frac{1}{2}}
		+ \frac{h}{2} \I^{-1} \left( \I\,W_{k+1} \times W_{k+1} + \tau(Q_{k+1}) \right)\;. \label{eq:LN_3}	&
\end{numcases}
Here we have introduced the Cayley map $\cay: \mathbb{R}^3 \to \SO{3}$:
\begin{equation} \label{eq:Cayley}
	\cay(x) =
		I + \frac{4}{4 + | x |}\hat{x} +
			\frac{2}{4 + | x |}{\hat{x}}^2 \;,~~~ x\in \Real{3} \;,
\end{equation}
where $I$ is the identity matrix.  The Cayley map is a second-order approximation
of the exponential map on $\SO{3}$.   There are other maps
one can use in place of the Cayley map in \eqref{eq:LN_2} (see, e.g., \cite[\S 5.4]{Bo2007}),
but the Cayley map is known to be very computationally efficient in practice.
This integrator is semi-explicit because \eqref{eq:LN_1}-\eqref{eq:LN_2} involve
explicit updates, and \eqref{eq:LN_3} is only implicit in the angular velocity and not in the torque.
Hence, the implicitness of the Lie-Newmark method is not severe.
It is also symmetric and second-order accurate.

\paragraph{Lie-Verlet}

The Lie-Verlet integrator was proposed in \cite{BoMa2008} and based on
the theory of discrete and continuous Euler-Poincar\'e systems \cite{MaPeSh1998,MaSc1993}.
The method is closely related to, but different from the RATTLE method for constrained
mechanical systems \cite{HaLuWa2006}.

Given  $(Q_k, W_k) \in \Real{3} \times \SO{3}$ and time-stepsize $h$, the Lie-Verlet algorithm
determines $(Q_{k+1}, W_{k+1})$ by the following iteration rule:
\begin{numcases} {\label{eq:LieVerlet}}
	W_{k+\frac{1}{2}} \! = \! W_k + \frac{h}{2} \I^{-1} \!\!\left[
		\I\,W_{k+\frac{1}{2}} \! \times\!  W_{k+\frac{1}{2}}
		- \frac{h}{2} \left( W_{k+\frac{1}{2}}^T \I\,W_{k+\frac{1}{2}}\right) W_{k+\frac{1}{2}}
		+ \tau(Q_k)\right] \;, & \label{eq:VLV_1} \\
	Q_{k+1} = Q_k\,\cay(h W_{k+\frac{1}{2}}) \;, & \label{eq:VLV_2} \\
	W_{k+1}\! = \! W_{k+\frac{1}{2}}\!\! + \!\frac{h}{2} \I^{-1}\!\!\left[
		 \I\,W_{k+\frac{1}{2}}\! \!\times \! W_{k+\frac{1}{2}} \!
		+ \! \frac{h}{2}\!\left( W_{k+\frac{1}{2}}^T \I\,W_{k+\frac{1}{2}} \right)\!W_{k+\frac{1}{2}} \!\!
		+ \tau(Q_{k+1}) \right]\!. &  \label{eq:VLV_3}
\end{numcases}
Similar to the Lie-Newmark method, this algorithm is symmetric, semi-explicit and second-order accurate.
In particular, the updates in \eqref{eq:VLV_2} and \eqref{eq:VLV_3} are explicit, and the implicitness
in \eqref{eq:VLV_1} does not involve the torque.  We emphasize the Lie-Verlet
integrator is variational and refer the reader to \cite{BoMa2008} for a proof of this result.


\section{Numerical Counterexample}

This section describes a numerical experiment
showing that the semi-explicit Lie-Newmark integrator
\eqref{eq:LieNewmark} exhibits systematic drift in
total energy.  Such drift implies that
the method is not a conjugate symplectic integrator for \eqref{eq:RigidBody},
and hence, is not variational.    The numerical counterexample we discuss
is strongly inspired by a numerical experiment reported in
\cite[\S 4.4]{FaHaPh2004}.
That paper shows systematic energy drift along an orbit of a
fourth-order accurate, implicit, and symmetric Lobatto IIIB  integrator
when applied to a spring pendulum with exterior forces.

\paragraph{Preliminaries}

Let $I \in \SO{3}$ be the identity matrix, and in terms of which,
define the function
\mbox{$\dist:\SO{3}\times\SO{3} \rightarrow \Real{}$} as
\begin{equation}
  \dist(Q_1,Q_2) \eqdef \sqrt{2\,\trace(I  - Q_1^T Q_2)} \,.
\end{equation}
Let $\|\cdot\|_F$ denote the Frobenius matrix norm.
We recall that $\|A\|_F := \sqrt{\trace(A^T A)}$
for $A \in \Real{n\times n}$.  It is straightforward to verify that 
$\dist(\cdot,\cdot)$ is a metric on $\SO{3}$ induced by the
Frobenius norm using the identity \[
 \| Q_2 - Q_1  \|^2_F = 2\,\trace(I - Q_1^T Q_2) \,.
\]

For the numerical experiment, consider
a single rigid body in a static field defined
by the potential energy function $U_{\alpha}:\SO{3} \rightarrow \Real{}$
given by:
\begin{equation} \label{eq:potential}
	U_{\alpha}(Q) = 
      \bigl( \dist(Q,I) - 1 \bigr)^2 
      - \frac{\alpha}{\dist(Q,Q_m)} \;.
\end{equation}
The first term in the right hand side of \eqref{eq:potential}
is a bounded potential which attains its minimum value
at $Q \in \SO{3}$ satisfying $\dist(Q,I) = 1$. 
The second term is an unbounded potential
that generates an attraction toward the
configuration $Q_m \in \SO{3}$.  
The parameter $\alpha$ is a tuning parameter.

For $\alpha = 0$, the potential energy $U_{\alpha=0}$ achieves
its minimum value on the two-dimensional surface 
\[
  S \eqdef \{Q \in \SO{3}: \dist(Q,I) = 1\} \;.
\]
This implies that the set
$S \times \{0\} \subset \SO{3}\times\Real{3}$
is a (locally) stable set in the sense of Lyapunov
for the dynamics of the rigid body.
One proves this fact using as Lyapunov function  
the energy  
$$
E(Q,W) = T(W) + U_{\alpha=0}(Q)$$
and noting that the set 
$\{(Q,W) \in \SO{3}\times\Real{3}\;|\; E(Q,W) \leq \bar E\}$
is compact for every $\bar E \geq 0$.

For $\alpha > 0$, the set $S$ gets perturbed
by the unbounded attractive potential.
On this perturbed energy landscape, 
the rigid body experiences an attraction 
toward the configuration $Q_m$.
Yet, if we place the attraction point $Q_m$ 
sufficiently far from the set $S$ and choose the tuning parameter $\alpha > 0$ 
sufficiently small, the set $S$ gets only slightly perturbed into a new set, 
that we label $S_\alpha$.  Furthermore, the set 
$S_\alpha \times \{0\} \subset \SO{3}\times\Real{3}$ is locally Lyapunov stable 
like the unperturbed set $S\times \{0\}$.

In summary, we can design $U_{\alpha}$ so that the true solution conserves a
compact energy function when initialized in a neighborhood of the set $S_\alpha \times \{0\}$.   
This conserved quantity implies that the true solution is confined to this neighborhood.  
This property will be preserved by the Lie-Verlet integrator since it is variational.   
Recall that a variational integrator is interpolated by a level set of a `modified' energy function nearby the
true energy \cite{BeGi1994, Re1999, HaLuWa2006}.   This implies that a trajectory of the
variational integrator is confined to a neighborhood of $S_\alpha \times \{0\}$ for the duration of the simulation.
Next we show that the Lie-Newmark integrator exhibits systematic drift in this 
energy function.

\paragraph{Numerical `Stress Test'}

Now we are in position to describe the numerical counterexample.
Select the inertia matrix to be \[
\I = \begin{bmatrix} 2 & 0 & 0 \\
                                  0 & 2 & 0 \\
                                  0 & 0 & 4 \end{bmatrix} \;,
\]
and the tuning parameter to be $\alpha = 0.3$.
Place the attraction point at \[
Q_m = \exp( v_m ) \;,~~~v_m = \frac{1}{\sqrt{2}} \begin{bmatrix} 2.5 \\ 0 \\  2.5 \end{bmatrix}  \;.
\]
Here $\exp: \Real{3} \to \SO{3}$ is the exponential map \[
\exp(x) = I + \frac{\sin( |x|)}{|x|} \hat{x} + \frac{1 - \cos(|x|)}{ |x|^2} \hat{x}^2 \;,~~~ x\in \Real{3} \;.
\]
The initial condition $(Q_0,W_0) \in \Real{3} \times \SO{3}$ is selected so
that $\dist(Q_0,I)$ is nearly one.  Specifically, select the initial configuration and angular velocity to be \[
Q_0 = \exp(v_0) \;,~~~v_0 = \begin{bmatrix} 0 \\ 0.7227 \\ 0 \end{bmatrix} \;,
~~~W_0 = \begin{bmatrix} 0 \\ 0 \\ 0.625 \end{bmatrix} \;.
\]

\begin{figure}[!h]%
\centering
	\subfigure[h = 0.125]{\label{fig:energyErrorA}
		\includegraphics[width=\figscale\columnwidth]{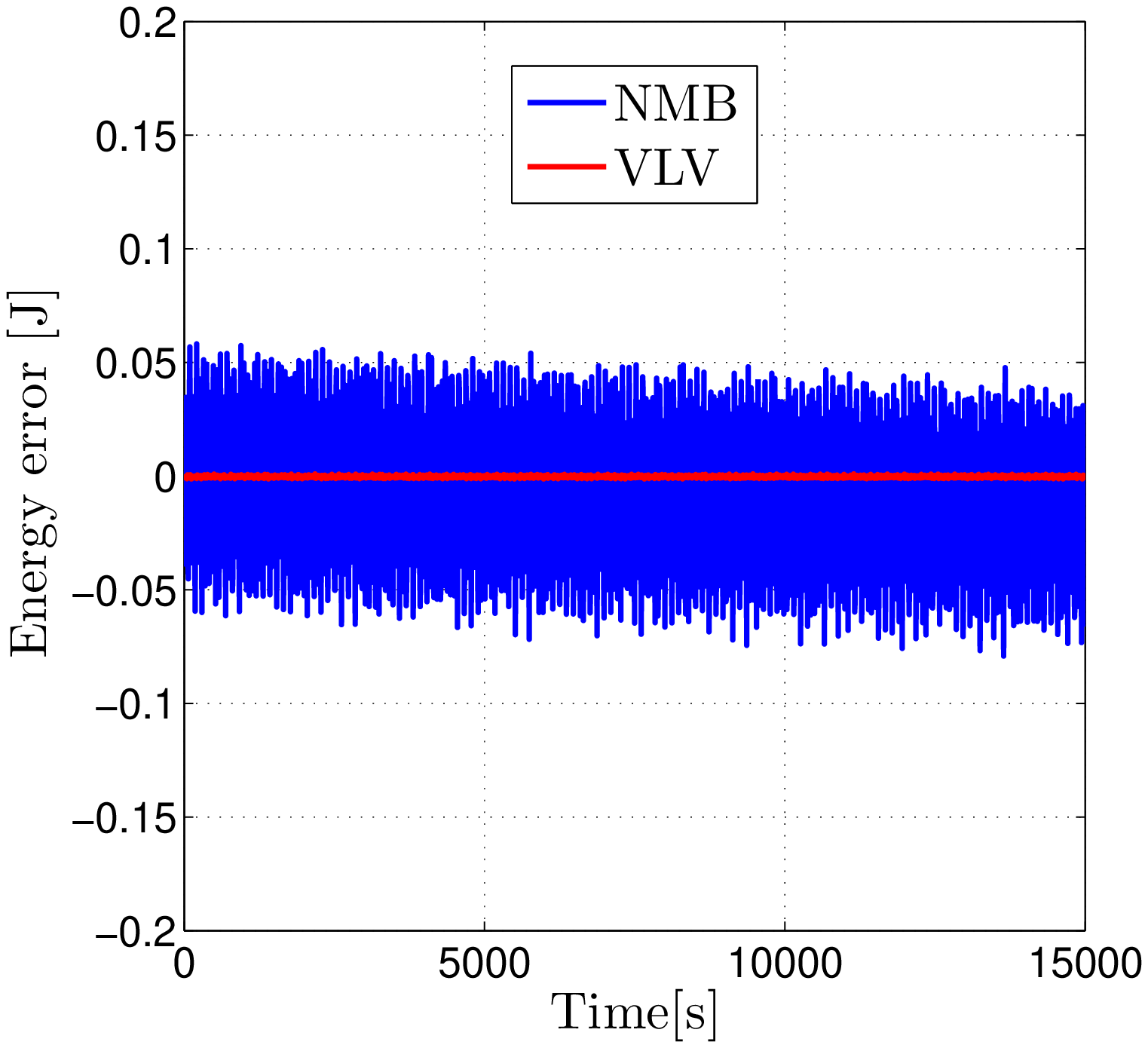} }
	\subfigure[h = 0.25]{\label{fig:energyErrorB}
		\includegraphics[width=\figscale\columnwidth]{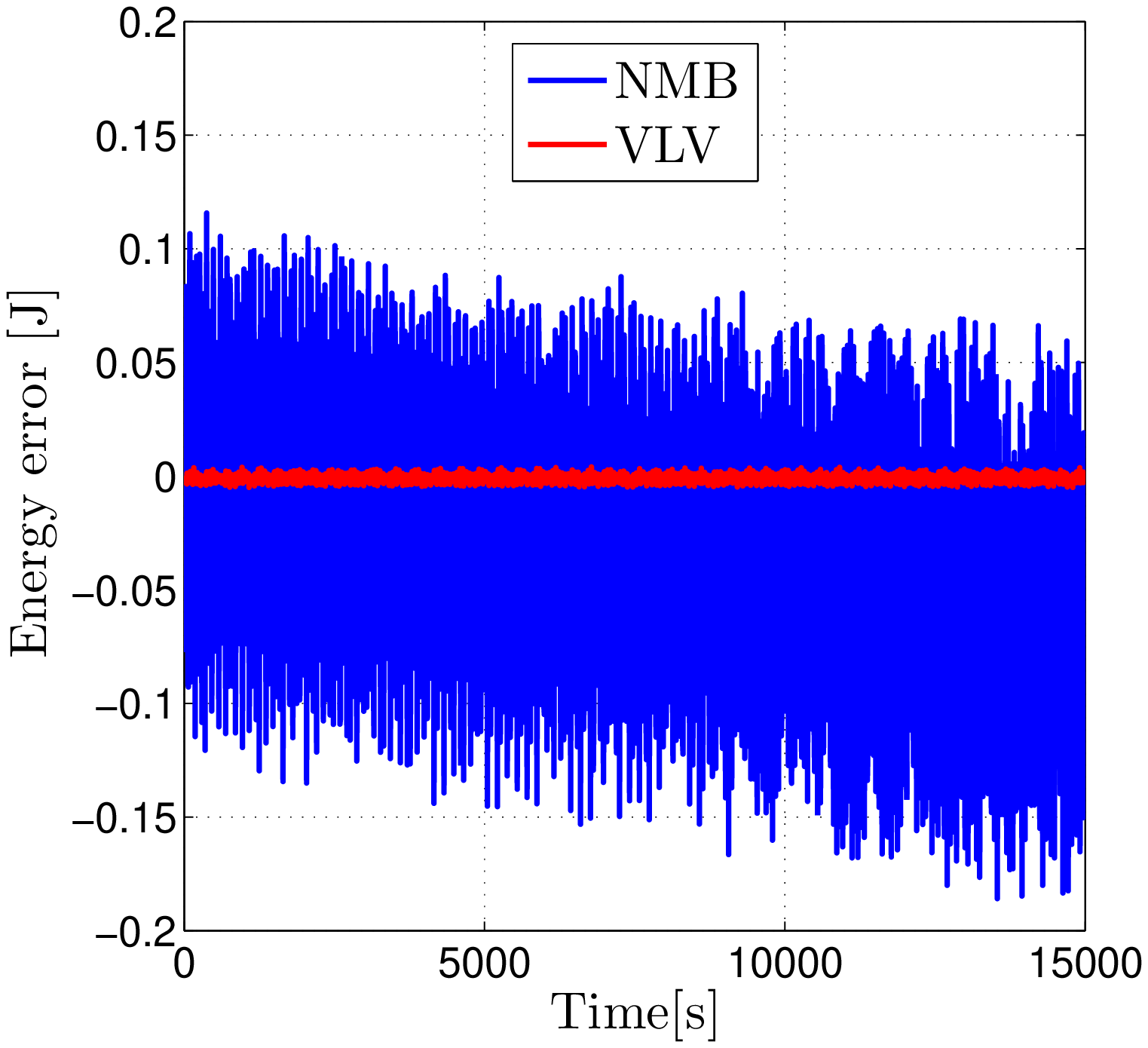}}
\caption{This figure shows the energy error of the Lie-Newmark (NMB) and velocity Lie-Verlet
(VLV) algorithms for the rigid body in the potential energy landscape defined by
\eqref{eq:potential} for  two different timesteps.  NMB exhibits a systematic energy drift.
On the other hand,  the energy error of VLV method remains bounded as predicted by theory.
The initial conditions and parameters used are provided in the text.
}%
\label{fig:energyError}%
\end{figure}

In the numerical experiment we test the two integrators,
Lie-Newmark (NMB) and velocity Lie-Verlet (VLV),
on a long time interval $[0, 15000]$.
The energy error obtained with a time-stepsize $h = 0.125$
is shown in Figure~\ref{fig:energyErrorA}.
The experiment was repeated with
a time-stepsize $h = 0.25$ and results are
reported in Figure~\ref{fig:energyErrorB}.
A systematic drift for the NMB scheme can be observed
in both cases.   The drift appears linear in
the time span $T$ and quadratic in the time-stepsize $h$.
We abbreviate this fact by saying the total energy
error behaves like $\mathcal{O}(T h^2)$.   No energy drift is
observed for the VLV scheme.
We have also tested
the Lie-Newmark method with the Cayley map in \eqref{eq:LN_2} replaced by the exponential
map.   Systematic energy drift is observed in that case too.

\begin{figure}[!h]
\centering{
	\includegraphics[width=\figscale\columnwidth]{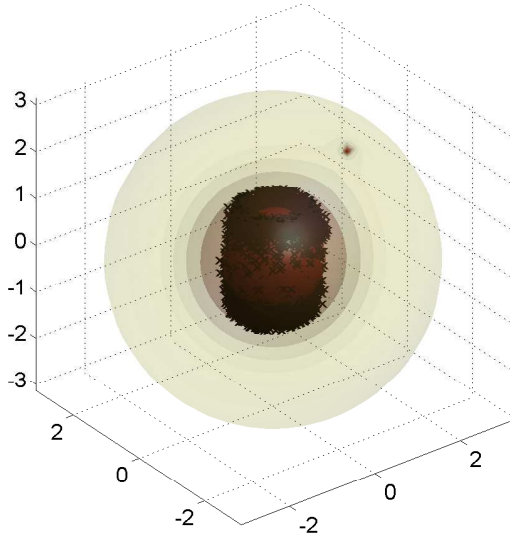}}
\caption{This figure shows a trajectory generated by the Lie-Newmark integrator operated 
at time-stepsize $h=0.25$ projected on the axis/angle representation of $\SO{3}$.  The 
initial conditions and parameters used are provided in the text.  The semi-transparent surfaces are
level sets of the potential energy \eqref{eq:potential}.  The darker shading corresponds 
to lower potential energy.  The dot in the figure corresponds to the attraction point 
$Q_m \in \SO{3}$ of the potential energy.
}%
\label{fig:trajectory}%
\end{figure}

 The trajectory generated by Lie-Newmark
for time-stepsize $h=0.25$ is shown in the \linebreak[4]axis/angle representation
of $\SO{3}$ in Figure~\ref{fig:trajectory}.  The semi-transparent
surfaces correspond to isosurfaces of the potential energy \eqref{eq:potential}.

\begin{figure}[!h]%
\centering
	\subfigure[Configuration.]{\label{fig:timeprecdiagrotation}
		\includegraphics[width=\figscale\columnwidth]{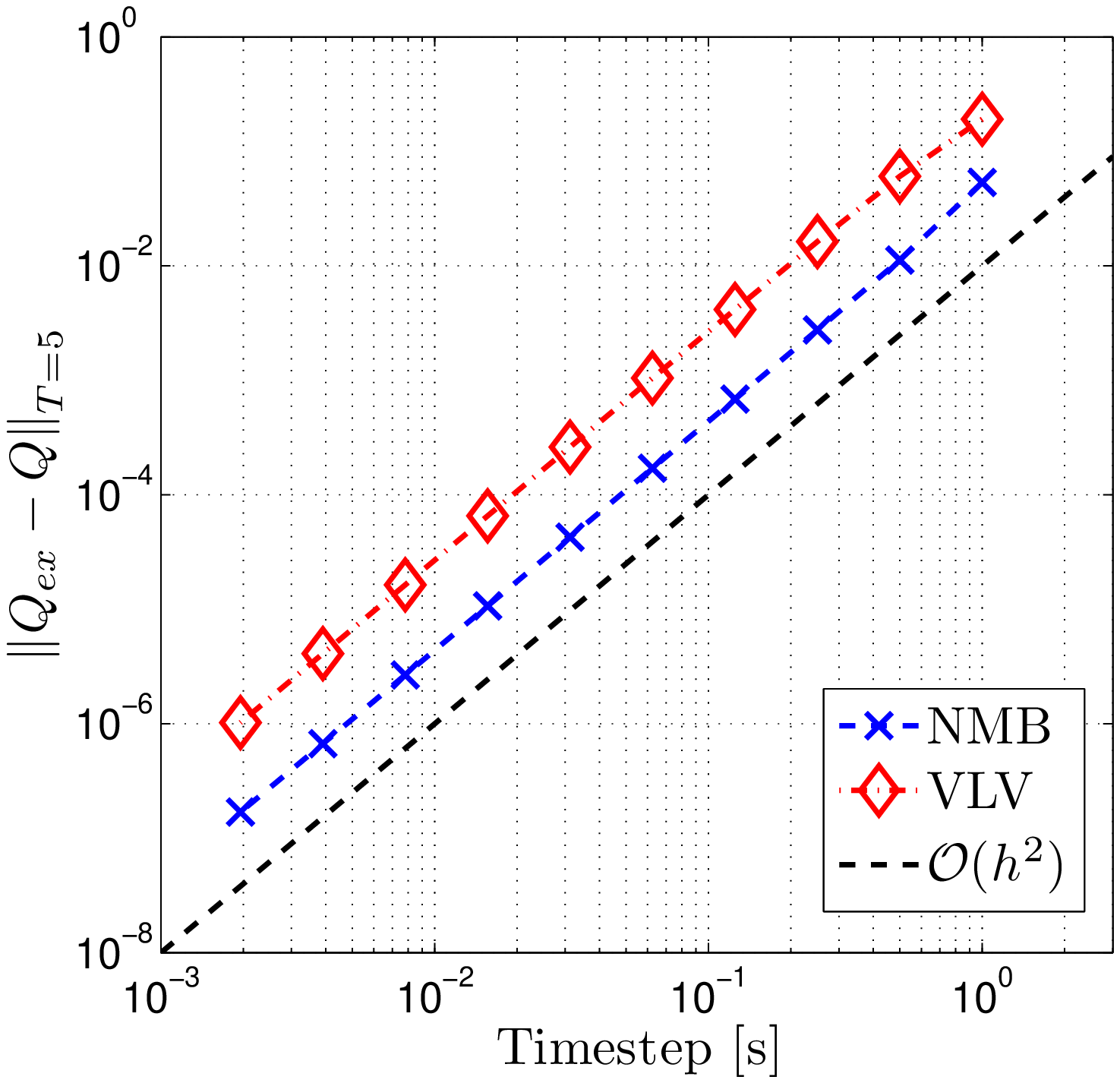} }
	\subfigure[Body angular velocity.]{\label{fig:timeprecdiagvelocity}
		\includegraphics[width=\figscale\columnwidth]{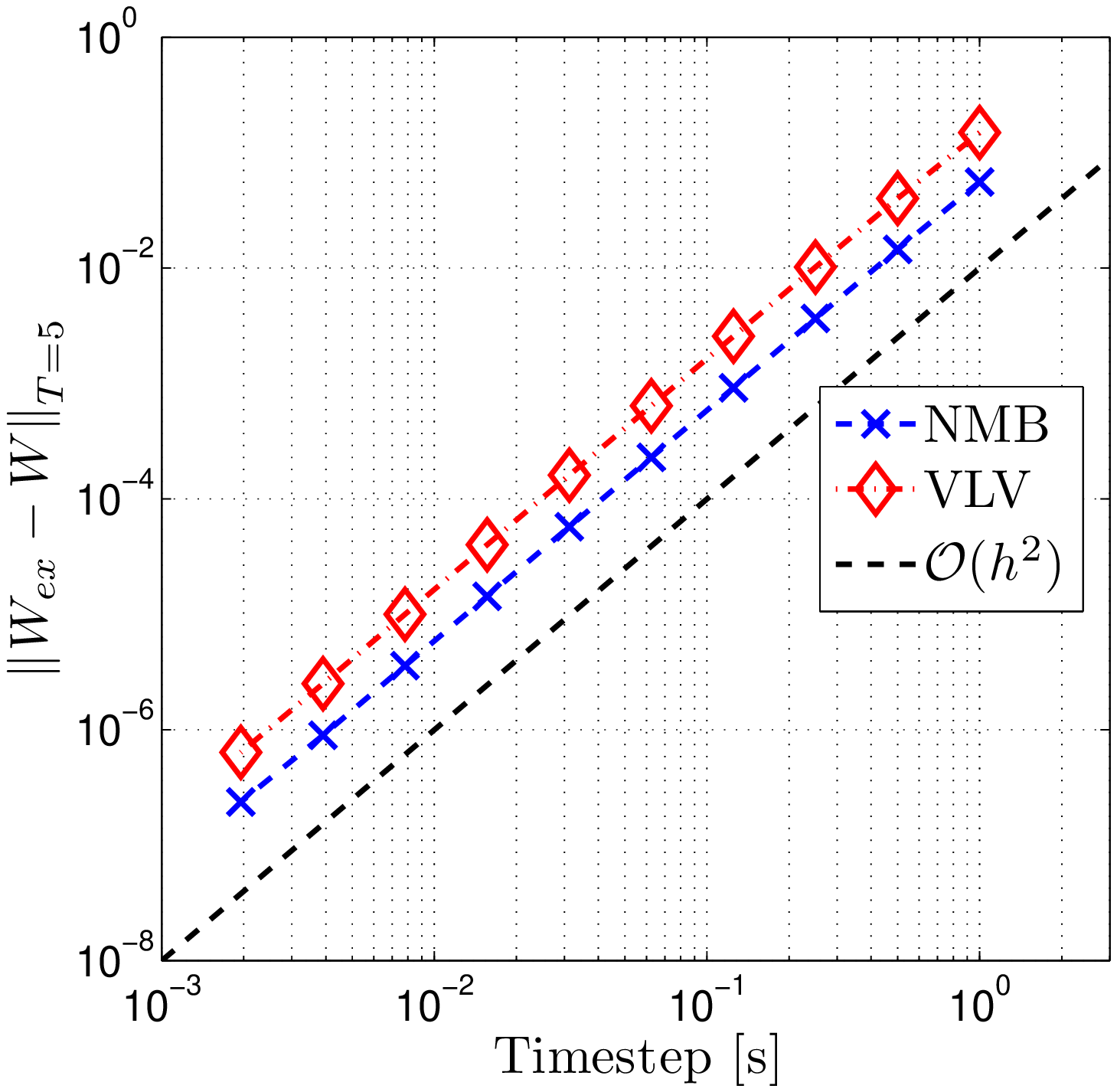}}
\caption{
This figure shows the global error of the Lie-Newmark (NMB) and velocity Lie-Verlet
(VLV) algorithms.  The global error is evaluated in configuration and body angular velocity
at a physical time of $T=5$ for a variety of time-stepsizes.   We use as a reference solution
an integration of \eqref{eq:RigidBody} using the \textsc{Matlab} function \texttt{ode45} with low
tolerance.  Observe that both integrators are second-order accurate.
}
\label{fig:timeprecdiag}%
\end{figure}

The time-precision diagrams,
shown in Figures~\ref{fig:timeprecdiagrotation} and \ref{fig:timeprecdiagvelocity}
confirm that NMB and VLV are second-order accurate.  Observe from the figures that
the slope of the two lines denoting the global error is $\mathcal{O}(h^2)$.
The diagrams have been generated by computing the global error
in the configuration and angular velocity evaluated at $T = 5$.
The simulations have been performed for a variety of time-stepsizes as indicated
in the figures.  The reference solution was computed using the
function \texttt{ode45} in \textsc{Matlab}, with an absolute tolerance $10^{-14}$
and relative tolerance $2\cdot 10^{-14}$.

\section{Conclusion}

The Lie-Newmark method was proposed as a generalization of the explicit Newmark algorithm to
Lie groups \cite{SiVu1988}.  However, unlike its counterpart on vector spaces, this paper shows that
the Lie-Newmark method does not possess excellent long-time behavior when applied to a rigid
body in a potential force field.  In particular, the paper presents a numerical experiment showing 
the Lie-Newmark integrator exhibits systematic energy drift with an $\mathcal{O}(T h^2)$ behavior.
In contrast, the Lie-Verlet  method, which is variational by construction, does not exhibit energy drift as
predicted by theory.   Since the two methods are semi-explicit, symmetric and computationally similar to implement,
we conclude that the Lie-Verlet method is better suited for long-time simulation of  rigid body-type systems.

As a final remark, we note that the numerical example appearing in this paper can serve as a test for symplecticity
of other rigid body integrators.  For instance, in Appendix \ref{app:LIEMIDEA} we show non-symplecticity
of the `explicit' Lie-Midpoint algorithm proposed in \cite{Kr2005}.


\section*{Acknowledgements}
We are extremely grateful to Jerry Marsden for suggesting this topic,
providing encouragement, excellent teaching, and many good ideas.  We
also wish to thank Petr Krysl and Melvin Leok for useful discussions.

NBR acknowledges the support
of the United States National Science
Foundation \linebreak[4]through NSF Fellowship \# DMS-0803095.
AS acknowledges the support of
the project \linebreak[4]DENO/FCT-PT (PTDC/EEA-ACR/67020/2006)
and the FCT-ISR/IST pluri-annual funding program.

\appendixpageoff
\appendixtitleon
\begin{appendices}


\section{Non-Symplecticity of Explicit Lie-Midpoint Algorithm}
\label{app:LIEMIDEA}

In this section the LIEMID[EA] algorithm is subjected to the numerical `stress test' described in \S 3.
This algorithm was thought to be a symplectic-momentum integrator in \cite{Kr2005}.
However, the test reveals that the integrator is not conjugate symplectic.

Given $(Q_k, W_k) \in \Real{3} \times \SO{3}$ and time-stepsize $h$, the LIEMID[EA] algorithm
determines $(Q_{k+1}, W_{k+1})$ by the following iteration rule:
\begin{numcases} {\label{eq:LIEMIDEA}}
	\Theta_{k+\frac{1}{2}} \! = \frac{h}{2} \I^{-1} \exp(-\frac{1}{2} \widehat{\Theta}_{k+\frac{1}{2}} ) (\I W_k + \frac{h}{2} \tau(Q_k))  \;, & \label{eq:LIEMIDEA_1} \\
	Q_{k+\frac{1}{2}} = Q_k \exp( \widehat{\Theta}_{k+\frac{1}{2}} )  \;, & \label{eq:LIEMIDEA_2} \\
	W_{k+\frac{1}{2}}\! = \! \I^{-1} \exp(- \widehat{\Theta}_{k+\frac{1}{2}}) (\I W_k + \frac{h}{2} \tau(Q_{k}) )  \;, &  \label{eq:LIEMIDEA_3} \\
			\Theta_{k+1} = \frac{h}{2} \I^{-1} \exp(-\frac{1}{2}  \widehat{\Theta}_{k+1}) (\I W_{k+\frac{1}{2}})  \;, &  \label{eq:LIEMIDEA_4} \\
			Q_{k+1}\! = Q_{k+\frac{1}{2}} \exp(  \widehat{\Theta}_{k+1})  \;, &  \label{eq:LIEMIDEA_5} \\
			W_{k+1}\! = \! \I^{-1} ( \exp(- \widehat{\Theta}_k) \I W_{k+\frac{1}{2}} + \frac{h}{2} \tau(Q_{k+1}) ) \! \;. &  \label{eq:LIEMIDEA_6}
\end{numcases}
This integrator is more implicit than the semi-explicit Lie-Newmark and Lie-Verlet methods
introduced in \S 2.  In particular, the updates \eqref{eq:LIEMIDEA_1} and \eqref{eq:LIEMIDEA_4}
are both implicit.  Hence, the algorithm involves two nonlinear solves per step, unlike the Lie-Newmark
and Lie-Verlet methods which involve one nonlinear solve per step.
The four remaining updates are explicit.  This algorithm was derived as a
composition of a half-step of a first-order Lie-midpoint method and its adjoint \cite{Kr2005}.   As an immediate
consequence, the method is symmetric and second-order accurate.
However, the numerical `stress test' shows that this integrator is not conjugate symplectic.

\begin{figure}[!h]%
\centering
	\subfigure[h = 0.125]{\label{fig:energyErrorC}
		\includegraphics[width=\figscale\columnwidth]{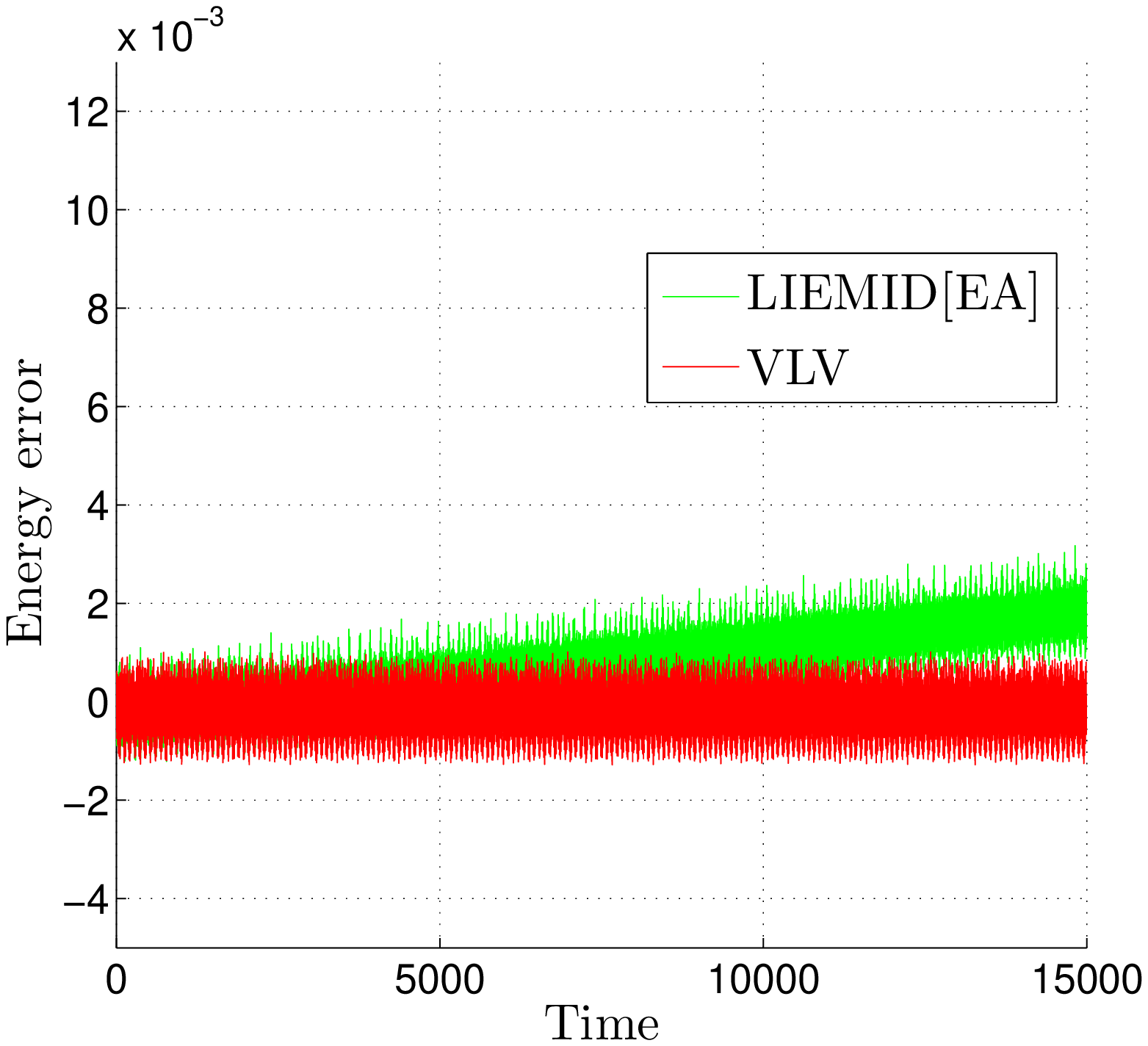} }
	\subfigure[h = 0.25]{\label{fig:energyErrorD}
		\includegraphics[width=\figscale\columnwidth]{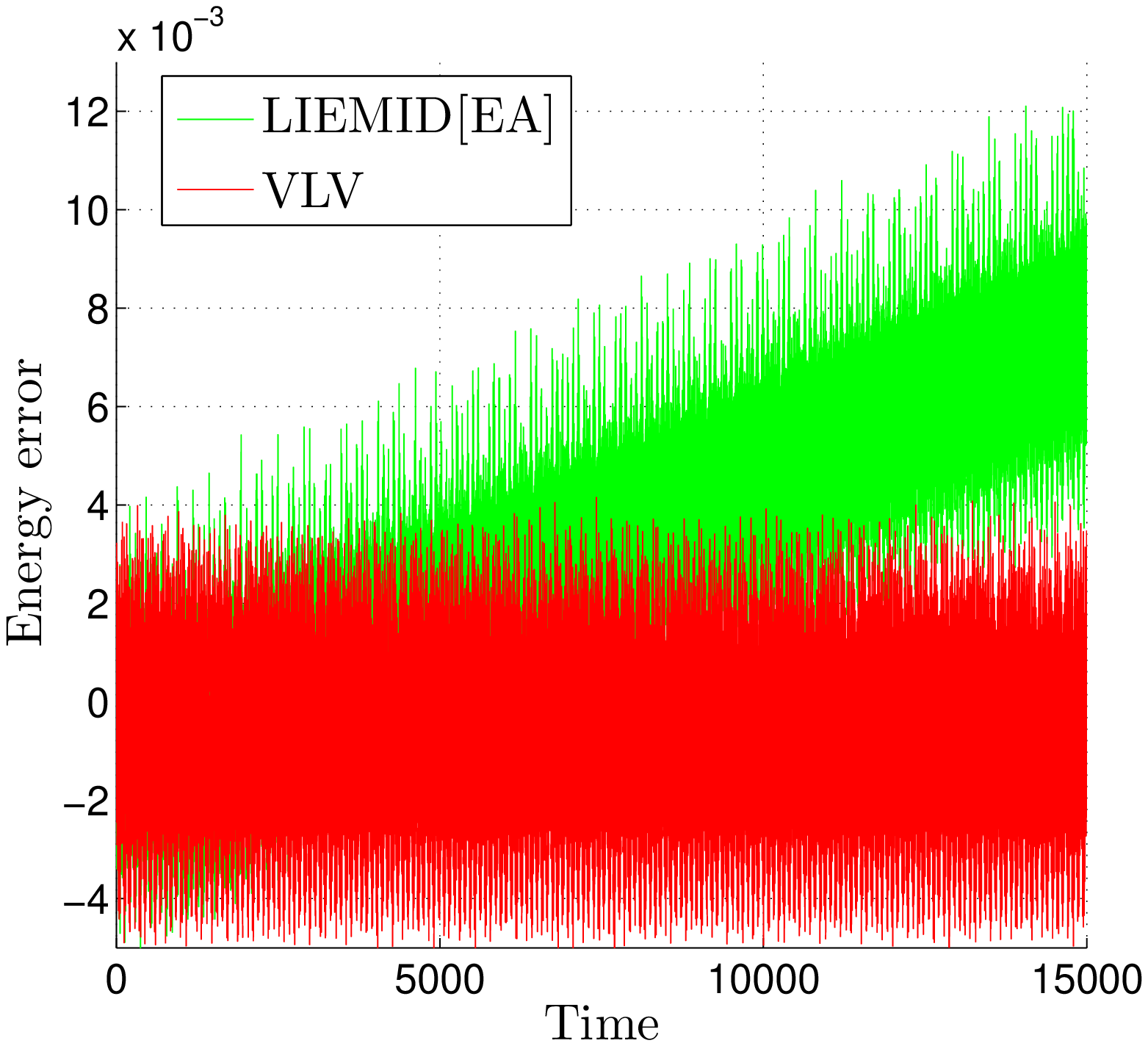}}
\caption{This figure shows the energy error of the Explicit Lie-Midpoint
(LIEMID[EA]) and velocity Lie-Verlet
(VLV) algorithms for the rigid body in the potential energy landscape defined by
\eqref{eq:potential} for two different timesteps.  LIEMID[EA] exhibits a systematic energy drift.
On the other hand,  the energy error of VLV method remains bounded as predicted by theory.
The initial conditions and parameters used are provided in the text.
}%
\label{fig:energyErrorLIEMID}%
\end{figure}

\begin{figure}[!h]%
\centering
	\subfigure[Configuration.]{\label{fig:timeprecdiagrotationB}
		\includegraphics[width=\figscale\columnwidth]{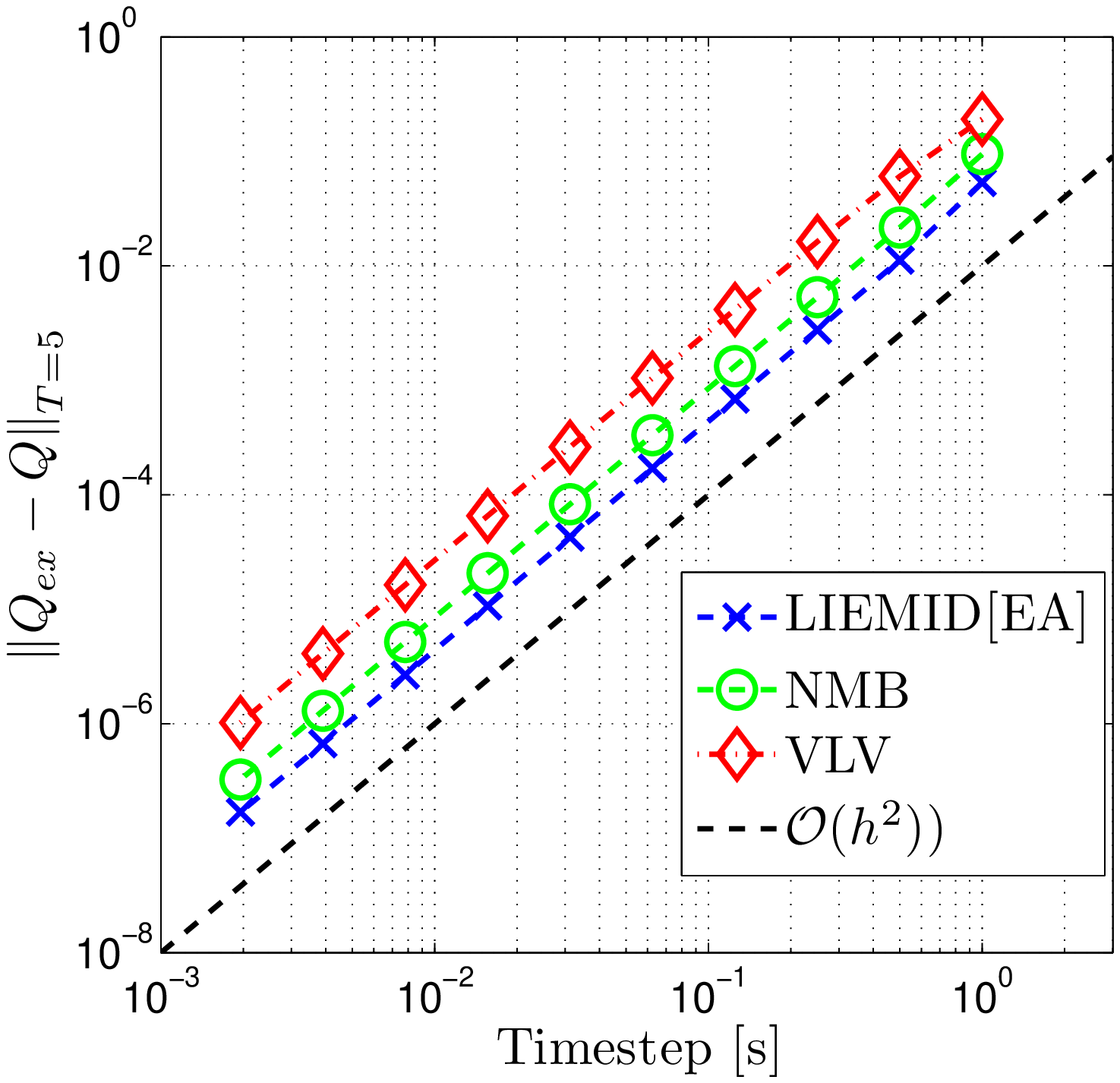} }
	\subfigure[Body angular velocity.]{\label{fig:timeprecdiagvelocityB}
		\includegraphics[width=\figscale\columnwidth]{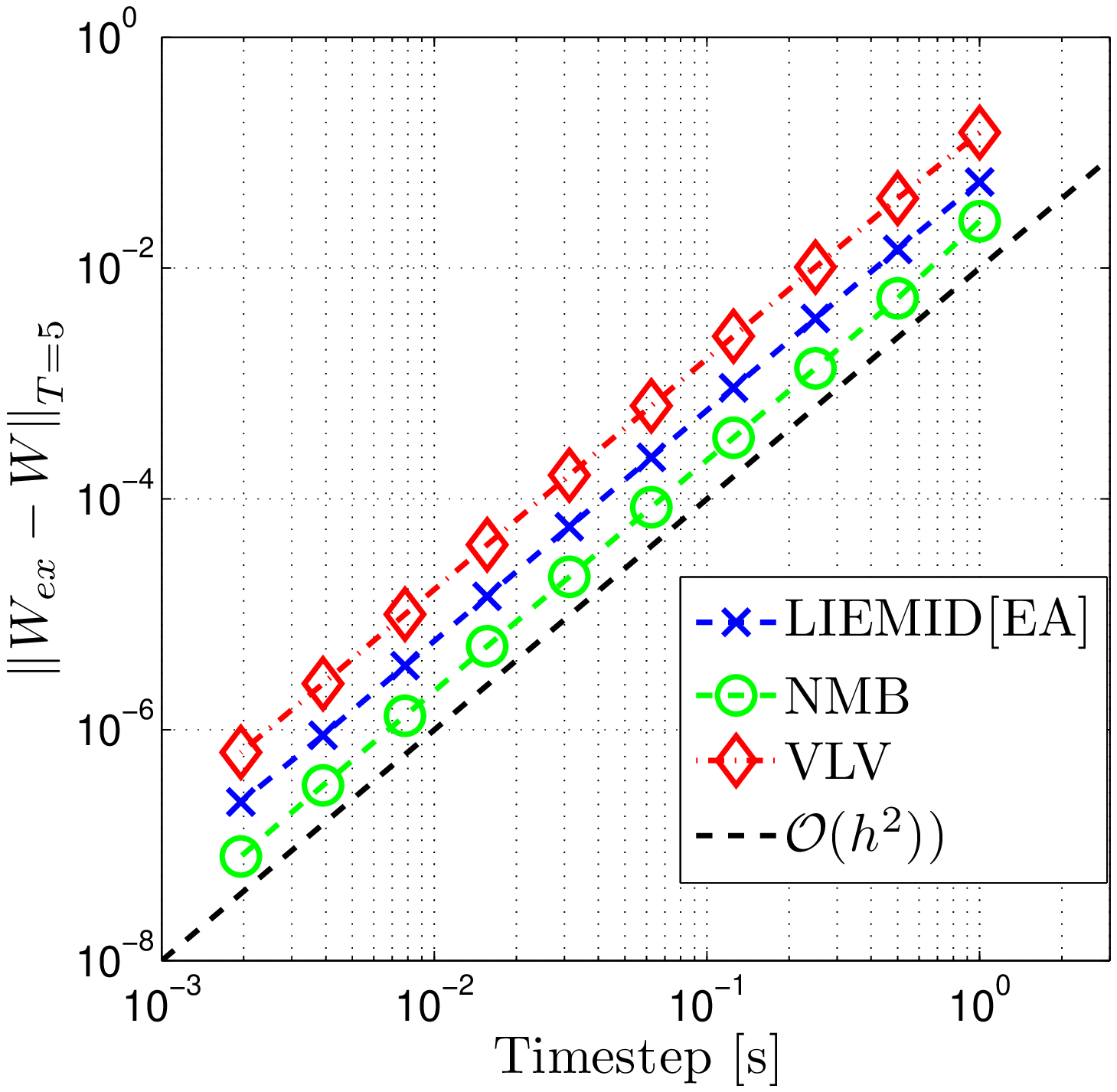}}
\caption{
This figure shows the global error of the Explicit Lie-Midpoint (LIEMID[EA]),
Lie-Newmark (NMB) and velocity Lie-Verlet (VLV) algorithms.
The global error is evaluated in configuration and body angular velocity
at a physical time of $T=5$ for a variety of time-stepsizes.   We use as a reference solution
an integration of \eqref{eq:RigidBody} using the \textsc{Matlab} function \texttt{ode45} with low
tolerance.  Observe that all the integrators are second-order accurate.
}
\label{fig:timeprecdiagALL}%
\end{figure}

The stress test described in \S 3 is carried out on LIEMID[EA].
The parameter values and initial conditions provided in \S 3 are used.
The time interval of integration is set to be $[0, 15000]$.
The energy error obtained on this time interval
with a time-stepsize $h = 0.125$ is shown in Figure~\ref{fig:energyErrorC}.
The experiment was repeated with
a time-stepsize $h = 0.25$ and results are
reported in Figure~\ref{fig:energyErrorD}.
A systematic drift for the LIEMID[EA] scheme can be observed
in both cases.   The drift appears linear in
the time span $T$ and quadratic in the time-stepsize $h$.
We note that the drift of LIEMID[EA] is smaller than the drift exhibited by NMB.
The time-precision diagrams,
shown in Figures~\ref{fig:timeprecdiagrotationB} and \ref{fig:timeprecdiagvelocityB}
confirm that LIEMID[EA] is second-order accurate.

\end{appendices}

\newpage

\bibliographystyle{Martin}
\bibliography{biblio}
\end{document}